\documentclass[12pt]{amsart}
\usepackage{etex}
\usepackage[english]{babel}
\usepackage{latexsym}
\usepackage{amsbsy}
\usepackage{amssymb}
\usepackage{amsfonts, amsmath}
\usepackage{mathtools}
\usepackage{csquotes}
\usepackage{pst-tree}
\usepackage{pst-math,pst-xkey}
\usepackage{multirow}
\usepackage{blkarray}
\usepackage{tabularx}
\usepackage{arydshln,leftidx,mathtools}
\usepackage{ifthen}
\usepackage{tikz}
\usepackage{tikz-cd}
\usepackage{float}
\usepackage{picins}
\usepackage{amsbsy}
\usetikzlibrary{matrix}
\usetikzlibrary{calc}
\usetikzlibrary{fit}
\usepackage{amsmath,amsfonts,amssymb}
\usepackage{mathrsfs,eurosym}
\usepackage{pst-plot,pst-eucl}
\usepackage{amsfonts}
\usepackage{alltt}
\usepackage{fontenc}
\usepackage{newlfont}
\usepackage{xlop}
\usepackage{graphicx}
\usepackage{yhmath}
\usepackage{mathdots}

\usepackage{MnSymbol}
\usepackage{amsthm}
\usepackage{url}
\usepackage{enumerate}
\usepackage{leqno}

\newtheorem{thm}{Theorem}[section]
\newtheorem{prop}[thm]{Proposition}

\newtheorem{ex}[thm]{Example}

\newtheorem{rmk}[thm]{Remark}

\theoremstyle{remark}

\numberwithin{equation}{section}
\theoremstyle{plain}

\DeclarePairedDelimiter\floor{\lfloor}{\rfloor}

\newcommand{\rdots}{\hspace{.2ex}\raisebox{1ex}{\rotatebox{-12}{$\ddots$}}}
\title[Powers of semicirculant and $r$-circulant matrices]{A closed-form expression for the kth power of semicirculant and $r$-circulant matrices}
\author{M. ~Mou\c{c}ouf}
\date{}
\begin{document}
\maketitle
\begin{center}
{\footnotesize Department of Mathematics, Faculty of Science, Chouaib Doukkali University, Morocco\\
Email: moucouf@hotmail.com}
\end{center}
\begin{abstract}
We derive a closed-form expression for the kth power of semicirculant matrices by using the determinant of certain matrices. As an application, a closed-form expression for the kth power of $r$-circulant matrices is also povided.
\end{abstract}
%
\section{Introduction}
%

The $n\times n$ $r$-circulant matrix $C_{n,r}$ over a unitary commutative ring $R$ is one having the
following form
\begin{equation}
C_{n,r}=\left(\begin{array}{cccccc}
c_{0}&c_{1}&c_{2}&\cdots&c_{n-2}&c_{n-1}\\
rc_{n-1}&c_{0}&c_{1}&\cdots&  c_{n-3}   &c_{n-2}\\
rc_{n-2}&rc_{n-1}&c_{0}&\cdots&c_{n-4}&c_{n-3}\\
\cdots&\cdots&\cdots&\cdots&\cdots  &\cdots\\
rc_{1}&rc_{2}&rc_{3}&\cdots&rc_{n-1}&c_{0}
\end{array}
\right),
\end{equation}
where $r\in R$ is a parameter. The $r$-circulant matrix $C_{n,r}$ is determined by $r$  and its first row elements $c_{0},\ldots,c_{n-1}$, thus we denote $C_{n,r}=\text{circ}_{n,r}(c_{0},\ldots,c_{n-1})$.
\\However, an infinite semicirculant matrix over $R$ is one having the following form (see e.g., Henrici~\cite{Hen} or Davis~\cite{Dav})
\begin{equation*}
A=[a_{0}, a_{1}, a_{2},\ldots]=
\begin{pmatrix}
a_{0}& a_{1}&a_{2}&a_{3}&\cdots\\
0&a_{0}&a_{1}&a_{2}&\cdots\\
0& 0&a_{0}&a_{1}& \cdots\\
0&0& 0&a_{0}& \cdots\\
\cdots&\cdots& \cdots&\cdots&\cdots \\
\end{pmatrix}.
\end{equation*}
In the work~\cite{Mou}, the general expression of the $k$th power $A^{k}=[a_{0}(k), a_{1}(k), a_{2}(k),\ldots]$ ($k\in \mathbb{N}$) of a semicirculant matrix
$A=[a_{0}, a_{1}, a_{2},\ldots]$ is presented. More precisely the sequence $\{a_{m}(k)\}_{m\geq0}$ is obtained using a recursive method. To do this, we have proved that for all $k\in \mathbb{N}$,
\begin{equation}\label{equ 13}
\begin{array}{ccc}
a_{i}(k) &=& L(A)(i,0)a_{0}^{k}
\binom{k}{0}+\cdots+L(A)(i,j) a_{0}^{k-j}
\binom{k}{j}+\cdots+L(A)(i,i) a_{0}^{k-i}
\binom{k}{i},
\end{array}
\end{equation}
where we adopt the convention that for any element $a\in R$ and any nonnegative integers $k\leq p$, $a^{k-p}\binom{k}{p}=\delta_{k,p}$.
\\The double sequence $\{L(A)(n,m)\}_{n,m}$ of elements of $R$ is defined by
\begin{eqnarray}\label{eq: 4001}
L(A)(n,m)=\displaystyle\sum_{\substack{\Delta(n,m)}}\binom{m}{k_{1},\ldots,k_{n}}a_{1}^{k_{1}}\cdots a_{n}^{k_{n}},
\end{eqnarray}
where $\Delta(n,m)$, $m\leqslant n$ are integers, is the solution set of the following system of equations
\begin{eqnarray*}
\begin{cases}
(k_{1},\ldots,k_{n})\in \mathbb{N}^{n}\\
k_{1}+\cdots+k_{n}=m\\
k_{1}+2k_{2}+ \cdots + nk_{n}=n.
\end{cases}
\end{eqnarray*}
We note that Formula \eqref{equ 13} expresses $a_{i}(k)$ in closed-form that is not easy to use. Here we provide an easy closed-form formula in terms of determinants. As a consequence, in the last section we give a closed-form formula for the kth power of $r$-circulant matrices and we describe a method for finding the solution set of the well known Diophantine equation
\begin{eqnarray*}
k_{1}+2k_{2}+ \cdots + nk_{n}=n.
\end{eqnarray*}
Throughout this paper $R$ will denote an arbitrary commutative ring with identity.
%
\section{Sequences $L(A)(n,m)$ and $u_{n}(A)$ corresponding to a semicirculant matrix $A$}
Let $A=[a_{0}, a_{1}, a_{2},\ldots]$ be a semicirculant matrix over $R$. Then The sequence $(L(A)(i,j))_{i,j}$ is independent on the coefficient $a_{0}$ and it is uniquely determined by the following recursive formula (see Lemma 2.2 of \cite{Mou}):
\begin{equation}\label{eq:SD2}
\left\{
\begin{array}{llllll}
L(A)(0,0)&=&1&\\
L(A)(i,0)&=&0& \quad \text{for } & i\neq 0\\
L(A)(i,j)&=&0& \quad \text{for } & i< j\\
L(A)(i,j+1)&=&a_{1}L(A)(i-1,j)+\cdots+a_{i-j}L(A)(j,j)
\end{array}
\right.
\end{equation}
Let $L(A)$ be the matrix $(L(A)(i,j))_{0\leqslant i,j}$. Then $L(A)$ is the following lower triangular matrix
$$ \begin{tikzpicture}[baseline=(math-axis),every right delimiter/.style={xshift=-3pt},every left delimiter/.style={xshift=3pt}]%
\matrix [matrix of math nodes] (matrix)
{
|(m11)| 1 &|(m12)| & |(m13)|& |(m14)|&|(m15)| &|(m16)|&|(m17)|&|(m18)|&|(m19)|\\
 |(m21)| 0  & |(m22)| a_{1} &|(m23)| &|(m24)| & |(m25)| &|(m26)| & |(m27)|&|(m28)| &|(m29)| \\
|(m31)| 0  & |(m32)| a_{2} & |(m33)| a_{1}^{2} & |(m34)|&|(m35)|     &|(m36)| &|(m37)| &|(m38)| &|(m39)| \\
|(m41)| 0  & |(m42)| a_{3} & |(m43)| 2a_{1}a_{2} & |(m44)| a_{1}^{3} & |(m45)| & |(m46)|   & |(m47)|& |(m48)|&|(m49)| \\
  |(m51)| 0  & |(m52)| a_{4} & |(m53)| 2a_{1}a_{3}+a_{2}^{2} & |(m54)| 3a_{1}^{2}a_{2}  & |(m55)| a_{1}^{4}    &|(m56)|   &|(m57)| &|(m58)|&|(m59)| \\
  |(m61)| 0  & |(m62)| a_{5} & |(m63)|2a_{1}a_{4}+2a_{2}a_{3} & |(m64)| 3a_{1}^{2}a_{3}+3a_{1}a_{2}^{2}  & |(m65)| 4a_{1}^{3}a_{2} & |(m66)|a_{1}^{5} & |(m67)|&|(m68)| &|(m69)| \\
  |(m71)|0  & |(m72)| a_{6} & |(m73)| 2a_{1}a_{5}+2a_{2}a_{4}+a_{3}^{2} & |(m74)|3a_{1}^{2}a_{4}+6a_{1}a_{2}a_{3}+a_{2}^{3} & |(m75)| 4a_{1}^{3}a_{3}+6a_{1}^{2}a_{2}^{2} & |(m76)| 5a_{1}^{4}a_{2} & |(m77)|a_{1}^{6} &|(m78)| &|(m79)|\\
  |(m81)|\vdots  &|(m82)|\vdots     & |(m83)|  \vdots & |(m84)|\vdots &|(m85)| \vdots
  &|(m86)| \vdots & |(m87)|\vdots       & |(m88)|\ddots & |(m89)|\\
};
\node[draw,dashed,inner sep=0.2pt,fit=(m33) (m63) (m33)] {};
\node[draw,dashed,inner sep=0.2pt,fit=(m74) (m74)] {};
\draw[ dashed ,->] (-2.6,0.2)--(-1.7,-1.3);
\coordinate (math-axis) at ($(matrix.center)+(0em,-0.25em)$);
\end{tikzpicture}$$
Each colomn of this matrix can be deduced from the precedent one.
For example $$\big<(a_{1}, a_{2}, a_{3}, a_{4}),(2a_{1}a_{4}+2a_{2}a_{3}, 2a_{1}a_{3}+a_{2}^{2}, 2a_{1}a_{2}, a_{1}^{2})\big>=3a_{1}^{2}a_{4}+6a_{1}a_{2}a_{3}+a_{2}^{3}.$$
Let $L_{i}(A)$ be the $i$th row of the matrix $L(A)$ and let $S_{i}(A)$ be the right shift of $L_{i}(A)$ by one position. Then formula~\eqref{eq:SD2} implies that
\begin{equation}\label{eq:SD3}
S_{n}(A)=a_{1}L_{n-1}(A)+\cdots+a_{n}L_{0}(A).
\end{equation}
We now consider the sequence
\begin{equation}\label{eq:SD4}
u_{n}(A)=\left\{
\begin{array}{llllll}
\sum_{i=1}^{n}L(A)(n,i)\\
u_{0}(A)=L(A)(0,0)=1
\end{array}\right.
\end{equation}
It is clear that $u_{n}$ is the sum of all elements of the ($n+1$)th row of $L(A)$. From formula~\eqref{eq:SD3} it follows that
\begin{equation}\label{eq:SD5}
\left\{
\begin{array}{llllll}
u_{n}(A)=a_{1}u_{n-1}(A)+a_{2}u_{n-2}(A)+\cdots+a_{n}u_{0}(A)\\
u_{0}(A)=1
\end{array}\right.
\end{equation}
Clearly, the sequence $u_{n}(A)$ is uniquely determined by the relation \eqref{eq:SD5}.
\\Let us now consider the following sets:
\\$\mathcal{SC}$ the set of all infinite semicirculant matrices with entries in the ring $R$,
\\$\mathcal{D}$ the subset of $\mathcal{SC}$ consisting of all diagonal matrices,
\\$\mathcal{U}$ the set of all sequences satisfying the recurrence relation~\eqref{eq:SD5} for some sequences $\{a_{n}\}_{n\geq1}$ of  the elements of $R$,
\\$\mathcal{V}$ the set of all double sequences satisfying the recurrence relation~\eqref{eq:SD2} for some sequences $\{a_{n}\}_{n\geq1}$ of the elements of $R$.\\
The following result shows the existence of a biunivoque correspondence between any two of the sets: $\mathcal{U}$, $\mathcal{V}$ and $\mathcal{SC}/\!\raisebox{-.65ex}{\ensuremath{\mathcal{D}}}$.
\begin{prop}\label{prop 18}
There is a commutative diagram of bijections
$$\begin{tikzcd}[row sep=2.5em]
 &\mathcal{SC}/\!\raisebox{-.65ex}{\ensuremath{\mathcal{D}}} \arrow{dr}{\varrho} \arrow{ld}[swap]{\psi}\\
 \mathcal{V}  \arrow{rr}{\varphi} && \mathcal{U}
\end{tikzcd}$$
where the maps $\rho$, $\psi$ and $\varphi$ are defined as follows
\begin{eqnarray*}
\varrho(\overline{A})&=&\{u_{n}(A)\}_{n\geq0}, \,\, \overline{A} \;\text{denotes the equivalence class of}\; A\in \mathcal{SC} \,\text{modulo}\, \mathcal{D}.\\
\psi(\overline{A})&=&\{L(A)(n,m)\}_{n,m},\\
\varphi(\{\phi(n,m)\}_{n,m})&=&\{\sum_{i=0}^{n}\phi(n,i)\}_{n\geq0}.
\end{eqnarray*}
These maps have respectively as inverses the following
\begin{eqnarray*}
\psi^{-1}(\{\phi(n,m)\}_{n,m})&=&\overline{[0,\phi(1,1),\ldots,\phi(n,1),\ldots]},\\
\varrho^{-1}(\{u_{n}\}_{n\geq0})&=&\overline{[0,a_{1},\ldots,a_{n},\ldots]},\\
\varphi^{-1}(\{u_{n}\}_{n\geq0})&=&\{\phi(n,m)\}_{n,m},
\end{eqnarray*}
where $\phi(n,m)=\displaystyle\sum_{\substack{\Delta(n,m)}}\binom{m}{k_{1},\ldots,k_{n}}a_{1}^{k_{1}}\cdots a_{n}^{k_{n}}$
and
$\{a_{n}\}_{n\geq1}$ is the sequence of elements of $R$ corresponding to $\{u_{n}\}_{n\geq0}$.
\end{prop}
\proof
Because sequences of $\mathcal{U}$ are uniquely determined by~\eqref{eq:SD5}, the map $\varrho$ is well defined. Let now $A=[a_{0},a_{1},a_{2},\ldots]$ and $B=[b_{0},b_{1},b_{2},\ldots]$ such that $u_{n}(A)=u_{n}(B),$ for all positive integer $n$. Then $u_{1}(A)=u_{1}(B)$ and hence $a_{1}=b_{1}$. Furthermore, the relation~\eqref{eq:SD5} implies that $a_{n}=a_{n}u_{0}(A)=u_{n}(A)-a_{1}u_{n-1}(A)-\cdots-a_{n-1}u_{1}(A)$. It then follows by an easy induction that $a_{n}=b_{n}$ for all positive integer $n$, i.e., $\overline{A}=\overline{B}$. Hence the map $\varrho$ is injective. To show that $\varrho$ is surjective, consider an element $(u_{n})_{n\in \mathbb{N}}$ of $\mathcal{U}$ and let $(a_{n})_{n\in \mathbb{N^{\ast}}}$ the corresponding sequence of $R$. It is clear that $\varrho(\overline{[0,a_{1},\ldots,a_{n},\ldots]})=\{u_{n}\}_{n\geq0}$, and then the map $\varrho$ is bijective.
\\The same argument as for $\varrho$ applies again to $\psi$.
\\Now let $\{\phi(n,m)\}_{n,m}\in \mathcal{V}$ and consider the matrix $A=[0,\phi(1,1),\ldots,\phi(n,1),\ldots]$. It is clear that $\varphi(\{\phi(n,m)\}_{n,m})=\{u_{n}(A)\}_{n\geq0}$ and then $\varphi$ is bijective.\\
The last assertion follows immediately from the fact that $\varphi\circ\psi= \rho$.
\endproof
\begin{prop}\label{prop 19} Let $a_{0},a_{1},a_{2},\ldots\in R$ and let $\delta(a_{0}, a_{1}, \ldots, a_{n})$ be the following determinant
\begin{equation}
\delta(a_{0}, a_{1}, \ldots, a_{n})=
\left|\begin{array}{cccccc}
a_{1}   & a_{2}        &\cdots  &\cdots   & \cdots  &a_{n}\\
a_{0}       &     \ddots    &\ddots  &         &        &\vdots      \\
0       &    \ddots     & \ddots &  \ddots &        &\vdots\\
\vdots  &    \ddots     & \ddots &\ddots   &  \ddots      &\vdots\\
\vdots  &               & \ddots &\ddots   & \ddots & a_{2}\\
0       & \cdots        & \cdots &   0     & a_{0}      &a_{1}
\end{array}\right|.
\end{equation}
Then we have
\begin{eqnarray}
\delta(a_{0}, a_{1}, \ldots, a_{n})&=&\delta(1, a_{1}, \ldots, a_{i}a_{0}^{i-1}, \ldots, a_{n}a_{0}^{n-1}) \label{eq:55}\\
                                &=&\delta(-1, a_{1}, \ldots, a_{i}(-a_{0})^{i-1}, \ldots,a_{n}(-a_{0})^{n-1}) \label{eq:56}
                                \end{eqnarray}
for all $n\in \mathbb{N^{\ast}}$.
\end{prop}
\proof The claim is trivially true for $a_{0}=0$. Suppose $a_{0}\neq 0$. For the proof of the first equality it suffices to multiply every column $C_{i}$ of the determinant $\delta(a_{0}, a_{1}, \ldots, a_{n})$ by $a_{0}^{i-1}$, and then multiply every row $R_{i}$ of the resulting determinant by $a_{0}^{-i+1}$.\\
To obtain the last equality, it suffices to multiply each of the even rows of the determinant $\delta(1, \ldots,a_{i}a_{0}^{i-1}, \ldots, a_{n}a_{0}^{n-1})$ by $-1$, and then multiply each of the even columns of the resulting determinant by $-1$.
\endproof
\begin{prop}\label{prop 20} Let $A=[a_{0},a_{1},a_{2},\ldots]$ be a semicirculant matrix over $R$, and let $\{u_{n}(A)\}_{n\geq0}$ be the sequence of elements of $R$ associated to $A$ given in~\eqref{eq:SD4}. Then we have
\begin{equation}\label{eq:54}
u_{n}(A)=\delta(-1, a_{1}, \ldots, a_{n})
\end{equation}
for all $n\in \mathbb{N^{\ast}}$.
\end{prop}
\proof
By~\eqref{eq:55} we know that $\delta(-1, a_{1},\ldots, a_{n})=\delta(1, a_{1}, \ldots, a_{i}(-1)^{i-1}, \ldots, a_{n}(-1)^{n-1})$. Assume that $n\geqslant1$, and let us expand the determinant $\delta(1, a_{1}, \ldots, a_{i}(-1)^{i-1}, \ldots, a_{n}(-1)^{n-1})$
along the first row. Then we get
\begin{eqnarray*}
\delta(-1, a_{1}, \ldots, a_{n})=\sum_{i=1}^{n-1}a_{i}\delta_{i},
\end{eqnarray*}
where $\delta_{i}$ is the cofactor associated with the entry $(-1)^{i-1}a_{i}$ of the matrix $\delta(1, a_{1}, \ldots, a_{i}(-1)^{i-1}, \ldots, a_{n}(-1)^{n-1})$.\\
It can be easily seen that the cofactor $\delta_{i}$ has the form
   \begin{eqnarray*}
   \delta_{i}=
   \left|
    \begin{array}{ccc;{6pt/6pt}ccclc}
     1& & \star & \multicolumn{4}{c}{\multirow{3}{0cm}{\huge$\star$}}&\\
      & \rdots& &&&&&\\
     0&&1&&&&&\\
     \hdashline[6pt/6pt]
     &&&a_{1}&-a_{2}&\cdots&a_{n-i}(-1)^{n+1-i}&\\
     &\multicolumn{1}{c}{\multirow{2}{0.5cm}{\huge$0$}}&&1&\ddots&&\vdots&\\
     &&& &\ddots&\ddots&\vdots&\\
     &&&0&&1&a_{1}&
    \end{array}
    \right|
   \end{eqnarray*}
   Thus for $i\geqslant1$ we have $\delta_{i}=\delta(-1,a_{1},\ldots,a_{n-i-1})$.
Hence
\begin{eqnarray*}
\delta(-1,a_{1},\ldots,a_{n})=\sum_{i=1}^{n-1}a_{i}\delta(-1,a_{1},\ldots,a_{n-i-1}).
 \end{eqnarray*}
It follows that the sequence $\delta(-1,a_{1},\ldots,a_{n})$ satisfies the recurrence relation~\eqref{eq:SD5}. But since $\delta(-1,a_{1})=a_{1}$, we have $\delta(-1,a_{1},\ldots,a_{n})=u_{n}$ for all $n\in \mathbb{N^{\ast}}$.
\endproof
The following proposition is useful for the proof of the main result.
\begin{prop}\label{prop 201} Let $X$ be an indeterminate over $R$ and let $\{a_{n}\}_{n\geq1}$ be a sequence of the elements of $R$. Let $\mathcal{X}_{n}(X)$ be the polynomial sequence defined by
\begin{equation}\label{eq:201}
\mathcal{X}_{n}(X)=\delta(-X,a_{1},\ldots,a_{n}).
\end{equation}
 Then one has
\begin{equation}\label{eq: 1201}
\mathcal{X}_{n}(X)=\displaystyle\sum_{i=1}^{n}L(A)(n,i)X^{n-i},
\end{equation}
where $A=[0,a_{1},a_{2},\ldots]$.
\end{prop}
\proof
Let $A=[0,a_{1},a_{2},\ldots]$ and consider the semicirculant matrix $B=[0,a_{1},a_{2}X,\ldots,a_{n}X^{n-1},\ldots]$ over $R[X]$. Then
\begin{eqnarray*}
L(B)(n,m)&=&\sum_{\substack{\Delta(n,m)}}\binom{m}{k_{1},\ldots,k_{n}}a_{1}^{k_{1}}(a_{2}X)^{k_{2}}\cdots (a_{n}X^{n-1})^{k_{n}}\\
&=&\sum_{\substack{\Delta(n,m)}}\binom{m}{k_{1},\ldots,k_{n}}a_{1}^{k_{1}}\cdots a_{n}^{k_{n}}X^{k_{2}+\cdots+(n-1)k_{n}},
\end{eqnarray*}
and since $k_{2}+\cdots+(n-1)k_{n}=(k_{1}+2k_{2}+\cdots+nk_{n})-(k_{1}+k_{2}+\cdots+k_{n})$ it follows that $k_{2}+\cdots+(n-1)k_{n}=n-m$, and $L(B)(n,m)=L(A)(n,m)X^{n-m}$. On the other hand, Formula~\eqref{eq:56} together with Formula~\eqref{eq:54} yield $u_{n}(B)=\mathcal{X}_{n}(X)$. As a consequence of Proposition~\ref{prop 18}, we have $\mathcal{X}_{n}(X)=\displaystyle\sum_{i=1}^{n}L(B)(n,i)=\sum_{i=1}^{n}L(A)(n,i)X^{n-i}$.
\endproof
Let $k$ and $n$ be nonnegative integers and let $[.]_{k}^{n}: R[X]\longrightarrow R[X]$ be the linear map defined by
\begin{eqnarray*}
[X^{i}]_{k}^{n}=X^{k-(n-i)}\binom{k}{n-i},
\end{eqnarray*}
where we use the convention that for $k\leq p$,
\begin{equation}\label{eq,,1}
X^{k-p}\binom{k}{p}=\delta_{k,p},
\end{equation}
the Kronecker delta. Then we have The following result which gives a closed-form expression of the kth power of semicirculant matrices.
\begin{thm}\label{Thmmm 1}
Let $\{a_{n}\}_{n\in \mathbb{N}}$ be a sequence of the elements of $R$. For all nonnegative integer $k$, the kth power of the semicirculant matrix $A=[a_{0},a_{1},a_{2},\ldots]$ are given as follows:
\begin{eqnarray}
A^{k}=[[\mathcal{X}_{0}]_{k}^{0}(a_{0}), [\mathcal{X}_{1}]_{k}^{1}(a_{0}), [\mathcal{X}_{2}]_{k}^{2}(a_{0}), \ldots]
\end{eqnarray}
where $\mathcal{X}_{n}$ is the polynomial $\delta(-X,a_{1},\ldots,a_{n})$ given in~\eqref{eq:201} and $\mathcal{X}_{0}=1$.
\end{thm}
\proof
Follows immediately from Proposition~\ref{prop 201} and Formula~\eqref{equ 13}.
\endproof
\begin{ex}
Consider a semicirculant matrix $A=[3, 2, 2, -3, 3]$. Let $k$ be any nonnegative integer.
We have
\begin{eqnarray*}
\mathcal{X}_{3}(X)= \begin{vmatrix} 2& 2& -3\\-X&2& 2\\0&-X&2\end{vmatrix}=-3X^{2}+8X+8.
\end{eqnarray*}
Then the $(1,4)$ entry of matrix $A^{k}$ is
$-3\times 3^{k-1}\binom{k}{1}+8\times 3^{k-2}\binom{k}{2}+8\times 3^{k-3}\binom{k}{3}$.
\end{ex}
\begin{ex}
Consider another semicirculant matrix $B=[0, 2, 1, 0,3]$. Let $k$ be any nonnegative integer.
We have
\begin{eqnarray*}
\mathcal{X}_{4}(X)= \begin{vmatrix} 2& 1&0&3\\-X&2&1&0\\0&-X&2&1\\0&0&-X&2\end{vmatrix}=3x^{3} + x^{2} + 12x + 16.
\end{eqnarray*}
Then the $(1,5)$ entry of matrix $B^{k}$ is $3\binom{k}{1}0^{k-1} + \binom{k}{2}0^{k-2} + 12\binom{k}{3}0^{k-3} + 16\binom{k}{4}0^{k-4}$.
\end{ex}
\section{Applications}
\subsection{A closed-form expression for the kth power of $C_{n,r}$}
\begin{thm}\label{Thm 10} Let $r, c_{0},\ldots,c_{n-1}\in R$ and let $k$ be any nonnegative integer. Consider the $r$-circulant matrix
\begin{eqnarray*}
C_{n,r}=\text{circ}_{n,r}(c_{0},\ldots,c_{n-1}).
\end{eqnarray*}
Then the pth strip of $C_{n,r}^{k}$ is
\begin{eqnarray*}
(C_{n,r}^{k})_{p}=\sum_{m\equiv p\;(\bmod\, n)}c_{0}^{k-m}[\mathcal{X}_{m}]_{k}^{m}(c_{0})r^{\floor*{\frac{m}{n}}},
\end{eqnarray*}
where
\begin{eqnarray*}
\mathcal{X}_{m}(X)=
\left\{
\begin{array}{lll}
\delta(-X,c_{1},\ldots,c_{m})& \text{if} &m\geq 1\\
1& \text{if}& m=0,
\end{array}
\right.
\end{eqnarray*}
and $\floor*{x}$ denotes the greatest integer less than or equal to $x$.
\end{thm}
\proof
Follows immediately from Theorem~\ref{Thmmm 1} and Theorem 4.1 of~\cite{Mou}.
\endproof
\begin{rmk}
The sequence $r^{\floor*{\frac{m}{n}}}$, appearing in the pth strip
\begin{eqnarray*}
(C_{n,r}^{k})_{p}=\sum_{m\equiv p}c_{m}(k)r^{\floor*{\frac{m}{n}}}
\end{eqnarray*}
of the $r$-circulant matrix $C_{n,r}^{k}$, is nothing but the geometric sequence with
common ratio $r$.
\end{rmk}
\begin{ex}
Let $C=circ_{5,2}(5, 4, 3, 2, 1)$ and let $[C]=[5, 4, 3, 2, 1,0,0,\ldots]$ be the associated infinite semicirculant matrix. Put
\begin{eqnarray*}
[C]^{k}=[c_{0}(k), c_{1}(k), c_{2}(k),\ldots].
\end{eqnarray*}
Using the method provided in Theorem~\ref{Thmmm 1}, we obtain
\begin{eqnarray*}
c_{s}(3)&=&0 \quad\text{for all}\quad s\geq (3\times4)+1=13;\\
\mathcal{X}_{12}(X)&=&x^{9} + \cdots, \quad c_{12}(3)=5^{9}\times5^{-9}\binom{3}{3}=1;\\
\mathcal{X}_{11}(X)&=&6x^{8} + \cdots, \quad c_{11}(3)=6\times 5^{0}\binom{3}{3}=6;\\
\mathcal{X}_{10}(X)&=&21X^{7} + \cdots,\quad c_{10}(3)=21\times 5^{0}\binom{3}{3}=21;\\
\mathcal{X}_{9}(X)&=&56X^{6} + \cdots,\quad c_{9}(3)=56\times5^{0}\binom{3}{3}=56;\\
\mathcal{X}_{8}(X)&=&x^{6} + 111x^{5} + \cdots,\quad c_{8}(3)=5^{1}\binom{3}{2}+111\times5^{0}\binom{3}{3}=126;\\
\mathcal{X}_{7}(X)&=&4X^{5} + 174X^{4} + \cdots,\quad c_{7}(3)=4\times5^{1}\binom{3}{2}+174\times5^{0}\binom{3}{3}=234;\\
\mathcal{X}_{6}(X)&=&10X^{4} + 219X^{3} + \cdots,\quad c_{6}(3)=10\times5^{1}\binom{3}{2}+219\times5^{0}\binom{3}{3}=369;\\
\mathcal{X}_{5}(X)&=&20X^{3} + 204X^{2} + \cdots,\quad c_{5}(3)=20\times5^{1}\binom{3}{2}+204\times5^{0}\binom{3}{3}=504;\\
\mathcal{X}_{4}(X)&=&x^{3} + 25x^{2} + 144x +\cdots,\quad c_{4}(3)=5^{2}\binom{3}{1}+25\times5^{1}\binom{3}{2}+144\times5^{0}\binom{3}{3}=594;\\
\mathcal{X}_{3}(X)&=&2X^{2} + 24X + 64,\quad c_{3}(3)=2\times5^{2}\binom{3}{1}+24\times5^{1}\binom{3}{2}+64\times5^{0}\binom{3}{3}=574;\\
\mathcal{X}_{2}(X)&=&3X + 16,\quad c_{2}(3)=3\times5^{2}\binom{3}{1}+16\times5^{1}\binom{3}{2})=465;\\
\mathcal{X}_{1}(X)&=&4,\quad c_{1}(3)=4\times5^{2}\binom{3}{1}=300;\\
c_{0}(3)=5^{3}&=&125.
\end{eqnarray*}
Hence \\$[C]^{3}=[125, 300, 465, 574, 594, 504, 369, 234, 126, 56, 21, 6, 1,0,0,\ldots]$.
\\Therefore,
\begin{eqnarray*}
C_{0}^{3}&=&125\times2^{0}+504\times2^{1}+21\times2^{2}=1217\\
C_{1}^{3}&=&300\times2^{0}+369\times2^{1}+6\times2^{2}=1062\\
C_{2}^{3}&=&465\times2^{0}+234\times2^{1}+2^{2}=937\\
C_{3}^{3}&=&574\times2^{0}+126\times2^{1}=826\\
C_{4}^{3}&=&594\times2^{0}+56\times2^{1}=706.
\end{eqnarray*}
Thus
\begin{eqnarray*}
C^{3}=circ_{5,2}(1217, 1062, 937, 826, 706).
\end{eqnarray*}
\end{ex}
\subsection{The calculation of $\Delta(n)$}
Consider the solution set $\Delta(n,p)$ of the following system of equations
\begin{eqnarray*}
\begin{cases}
(k_{1},\ldots,k_{n})\in \mathbb{N}^{n}\\
k_{1}+\cdots+k_{n}=p\\
k_{1}+2k_{2}+ \cdots + nk_{n}=n.
\end{cases}
\end{eqnarray*}
Here, $p$ and $n$ are positive integers such that $p\leqslant n$.
\\Let $\Delta(n)$ be the solution set of the diophantine equation
\begin{eqnarray*}
k_{1}+2k_{2}+ \cdots + nk_{n}=n.
\end{eqnarray*}
We have clearly
\begin{eqnarray*}
\Delta(n)=\bigcup_{1\leq i\leq n} \Delta(n,i).
\end{eqnarray*}
Let $X,x_{1},\ldots,x_{n}$ be independent indeterminates over $R$. From \eqref{eq: 1201} we have
\begin{eqnarray*}
\mathcal{X}_{n}(X)=\displaystyle\sum_{i=1}^{n}L(A)(n,i)X^{n-i},
\end{eqnarray*}
and from~\eqref{eq: 4001} we have
\begin{eqnarray*}
L(A)(n,i)=\displaystyle\sum_{\substack{\Delta(n,i)}}\binom{i}{k_{1},\ldots,k_{n}}x_{1}^{k_{1}}\cdots x_{n}^{k_{n}},
\end{eqnarray*}
where $A=[0,x_{1},\ldots,x_{n}]$.
\\We observe that we can determine the set $\Delta(n)$ by computing the polynomial $\mathcal{X}_{n}(X)$, and then Comparing the polynomial $L(A)(n,i), 1\leq i\leq n,$ and the term of degree $n-i$ of $\mathcal{X}_{n}(X)$.
\begin{ex}
We have $\mathcal{X}_{4}(X)=x_{4}X^{3}+(2x_{1}x_{3}+x_{2}^{2})X^{2}+3x_{1}^{2}x_{2}X+x_{1}^{4}$. Then
\begin{eqnarray*}
\Delta(4,1)&=&\{(0,0,0,1)\},\\
\Delta(4,2)&=&\{(1,0,1,0),(0,2,0,0)\},\\
\Delta(4,3)&=&\{(2,1,0,0)\},\\
\Delta(4,4)&=&\{(4,0,0,0)\}
\end{eqnarray*}
\end{ex}
\begin{ex}
We have $\mathcal{X}_{6}(X)=x_{6} X^{5} + (2x_{1}x_{5} + 2x_{2}x_{4} + x_{3}^{2})X^4 + (3x_{1}^{2}x_{4} + 6x_{1}x_{2}x_{3} + x_{2}^{3})X^{3}+ (4x_{1}^{3}x_{3} + 6x_{1}^{2}x_{2}^{2})X^{2} + 5x_{1}^{4}x_{2}X + x_{1}^{6}$. Then
\begin{eqnarray*}
\Delta(6,1)&=&\{(0,0,0,0,0,1)\},\\
\Delta(6,2)&=&\{(1,0,0,0,1,0),(0,1,0,1,0,0),(0,0,2,0,0,0)\},\\
\Delta(6,3)&=&\{(2,0,0,1,0,0),(1,1,1,0,0,0),(0,3,0,0,0,0)\},\\
\Delta(6,4)&=&\{(3,0,1,0,0,0),(2,2,0,0,0,0)\},\\
\Delta(6,5)&=&\{(4,1,0,0,0,0)\},\\
 \Delta(6,6)&=&\{(6,0,0,0,0,0)\}
\end{eqnarray*}
\end{ex}


\begin{thebibliography}{99}
\bibitem{Mou} M. Mou\c{c}ouf, \emph{Arbitrary positive power of semicirculant and $r$-circulant matrices,} 	arXiv:2006.15048 [math.RA].
\bibitem{Dav} P. Davis, \emph{Circulant matrices,} Ams Chelsea Publishing, Providence, Rhode Island, 2012.
\bibitem{Hen} P. Henrici, \emph{Applied and Computational Complex analysis,} volume 1, John Wiley, New York, 1974.
\end{thebibliography}
\end{document}